\documentclass{amsart}
\usepackage{amsmath,amsthm}
\usepackage{graphics}
\usepackage{color}
\usepackage{epsfig}

\graphicspath{{./Figures/}}

\theoremstyle{plain}
\newtheorem{thm}{Theorem}[section]
\newtheorem{cor}[thm]{Corollary}

\newtheorem{prop}[thm]{Proposition}
\newtheorem{lemma}[thm]{Lemma}

\newtheorem*{question}{Question}

\theoremstyle{definition}
\newtheorem{remark}[thm]{Remark}

\newcommand{\B}{\ensuremath{{\mathfrak b}}}
\newcommand{\comment}[1]{}

\newcommand{\bdry}{\ensuremath{\partial}}

\newcommand{\Z}{\ensuremath{\mathbb{Z}}}

\begin{document}
\title{Counting genus one fibered knots in lens spaces}

\author{Kenneth L. Baker}
\address{Department of Mathematics, University of Georgia \\ Athens, Georgia 30602}
\email{kb@math.uga.edu}

\thanks{This work was partially supported by a VIGRE postdoc under NSF grant number DMS-0089927 to the University of Georgia at Athens.}



\begin{abstract}
The braid axis of a closed $3$--braid lifts to a genus one fibered knot in the double cover of $S^3$ branched over the closed braid.  Every (null homologous) genus one fibered knot in a $3$--manifold may be obtained in this way.  Using this perspective we answer a question of Morimoto about the number of genus one fibered knots in lens spaces.  %We determine the number of genus one fibered knots in any given lens space and show that this number is $3$ in the case of the lens space $L(4,1)$ and at most $2$ otherwise.
We determine the number of genus one fibered knots up to homeomorphism in any given lens space.  This number is $3$ in the case of the lens space $L(4,1)$, $2$ for the lens spaces $L(m,1)$ with $m>0$, and at most $1$ otherwise.
\end{abstract}

\maketitle

\section{Introduction}
Let $M$ be a $3$--manifold.  We say a knot $K$ in $M$ is a {\em genus one fibered knot}, GOF-knot for short, if $M-N(K)$ is a once-punctured torus bundle over the circle whose monodromy is the identity on the boundary of the fiber and $K$ is ambient isotopic in $M$ to the boundary of a fiber.  In particular we will always consider a GOF-knot to be null homologous. 

As begun by Burde and Zieschang, Gonz\'alez-Acu\~na shows that the trefoil (and its mirror) and the figure eight knot are the only GOF-knots in $S^3$, \cite{bz:nkuf} and \cite{ga:dcok} respectively.  Morimoto shows that each lens space $L(m,1)$ contains at least two GOF-knots if $m > 0$ and exactly two if $m \in \{1,2,3,5,19\}$, $L(4,1)$ contains exactly three GOF-knots, each $L(0,1)$, $L(5,2)$, and $L(19,3)$ contains exactly one GOF-knot, and $L(19,2)$, $L(19,4)$, and $L(19,7)$ contain no GOF-knots, \cite{morimoto:gofkils}.  Morimoto then asks the following question.
\begin{question}[\cite{morimoto:gofkils}]
Are the numbers of GOF-knots in all lens spaces bounded?
\end{question}
In this article we use double branched covers of two-bridge links represented as closed $3$--braids to address this question.  

\begin{cor}\label{cor:boundgofkils}
Each lens space $L(m,1)$ with $m>0$ contains exactly two GOF-knots except $L(4,1)$ which contains three.  All other lens spaces contain at most one GOF-knot.

\end{cor}

This corollary is simplified version of Corollary~\ref{cor:countgofkils} where using Theorem~\ref{thm:ratlreps} we count the number of GOF-knots in any given lens space.

Throughout this article we will be considering links and $3$--manifolds up to homeomorphism.  Therefore, for instance, we regard the right-handed trefoil and left-handed trefoil in $S^3$ as equivalent.  In general, except for where noted, we consider links to be equivalent to their mirror and do not distinguish the orientations on a $3$--manifold.

If $h$ is a homeomorphism between $3$--manifolds $M$ and $M'$ such that $h(K)=K'$ for knots $K \subset M$ and $K' \subset M'$, then we say the pairs $(M,K)$ and $(M',K')$ are equivalent or simply that the knots $K$ and $K'$ are equivalent.  If $h$ is a homeomorphism of $S^3$ such that $h(L) = L'$ and $h(A)=A'$ for links $L$ and $L'$ with axes $A$ and $A'$ giving closed braid representations of $L$ and $L'$ respectively, then we say the pairs $(L,A)$ and $(L',A')$ are equivalent.  We further say $A$ and $A'$ are equivalent axes for $L$.

 Let $\omega$ be the braid word of a braid whose closure is the link $L$ with braid axis $A$.  Observe that the absolute value of the exponent sum of $\omega$ is an invariant of the equivalence class of the pair $(L,A)$.

We refer the reader to \cite{bz:knots} for background regarding fibered knots, braids,  two-bridge links, lens spaces, and double coverings of $S^3$ branched over a link.

%%%%%%%%%%%%%%%%%%
\section{GOF-knots via double branched covers of closed $3$--braids}

Each orientable once-punctured torus bundle is the double cover of a solid torus branched over a closed $3$--braid.  Moreover a Dehn filling of such a once-punctured torus bundle along a slope that intersects each fiber once is the double cover of $S^3$ branched over a closed $3$--braid.  This may be seen as follows.  The once-punctured torus $T$ admits an involution $\tau$ with three fixed points.  Any orientation preserving, boundary fixing homeomorphism of $T$ is isotopic rel--$\bdry$ to one invariant under $\tau$.  The involution $\tau$ then extends across each fiber of the mapping torus of such a homeomorphism.  The involution further extends across the solid torus of the Dehn filling described above (cf.\ Sections 4 and 5 \cite{mr:gfm}).

A GOF-knot in a $3$--manifold $M$ is then the lift of the braid axis of some closed $3$--braid in $S^3$ where $M$ is the double cover of $S^3$ branched over the closed $3$--braid.  This also becomes evident by considering the involution on the genus 2 Heegaard splitting of $M$ induced by the GOF-knot (cf.\ Section 5, \cite{birman-hilden:hsobcoS3}).

We codify this observation in the following theorem.

\begin{thm}\label{thm:generalgofks}

Up to equivalence\comment{homeomorphism}, the pairs $(M,K)$ of a GOF-knot $K$ in a $3$--manifold $M$ are in one-to-one correspondence with the pairs $(L,A)$ of a link $L$ in $S^3$ and a braid axis $A$ giving a closed $3$--braid representation of $L$.  The pair $(M,K)$ corresponds to the pair $(L,A)$ if and only if $M$ is the double cover of $S^3$ branched over $L$ and $K$ is the lift of $A$.
\end{thm}

\begin{remark}\label{remark:axissurgery}
A meridian of the braid axis is a longitudinal curve on the solid torus containing the closed $3$--braid and lifts to two meridians of the GOF-knot in the double cover.  More generally, let $V$ be the solid torus neighborhood of the braid axis and $\tilde{V}$ be the solid torus neighborhood of the GOF-knot that is the lift of $V$.  The meridian of the solid torus $S^3-\int(V)$ is the longitude of $V$ and lifts to the longitude of $\tilde{V}$.  Since $\tilde{V}$ double covers $V$, simple closed curves of slope $p/q$ on $\bdry V$ lift to curves of slope $2p/q$ on $\bdry \tilde{V}$.  Assuming $p$ and $q$ are coprime, if $q$ is even the slope $2p/q$ is to be interpreted as two parallel curves of slope $p/(q/2)$.  It follows that $1/n$ surgery on GOF-knot corresponds to inserting $2n$ full twists (right-handed if $n<0$, left-handed if $n>0$) into the $3$--braid.
%Slopes are taken with respect to standard meridian-longitude pairs.
\end{remark}

\begin{lemma}\label{lem:boundaxes}
An unoriented link $L$ has at most four equivalence classes of braid axes that represent $L$ as a closed $3$--braid.
\end{lemma}

\begin{proof}
By Lemma~\ref{lem:distinctbraids} below, up to reversal each orientation of a link $L$ admits at most one equivalence class of braid axes representing the oriented link as a closed $3$--braid except when $L$ is a type $(2,k)$ torus link with $k>0$.  When $L$ is oriented as a type $(2,k)$ torus link with $k>0$ it has two distinct equivalence class of braid axes representing it as a closed $3$--braid.

Since a $(2,k)$ torus link has at most two components, it has at most two distinct orientations up to reversal.  Therefore it may have at most four distinct equivalence classes of braid axes representing it as a closed $3$--braid. 

A link $L$ which may be represented as a closed $3$--braid has at most three components, and so has at most four distinct orientations up to reversal.  Assuming $L$ is not a $(2,k)$ torus link, each orientation up to reversal admits at most one equivalence class of braid axes representing it as a closed $3$--braid.  Therefore $L$ may have at most four equivalence classes of braid axes representing it as a closed $3$--braid. 
\end{proof}

For any given $3$--manifold $M$, there may be several different links in $S^3$ with $M$ as their double branched covers (\cite{bgam:hsop3manu}, \cite{viro:ltsbcab}, \cite{bedient:dbcapk}, among others).  In particular, there may be several different links with representations as closed $3$--braids that have $M$ as their double branched covers. 

The lens space $L(\alpha, \beta)$ is the double cover of $S^3$ branched along the (unoriented) two-bridge link $\B(\alpha, \beta)$; see, e.g.,  \cite{bz:knots}.  By classifying involutions on lens spaces Hodgson and Rubinstein show that if the lens space $L(\alpha, \beta)$ is the double cover of $S^3$ branched over a link $L$ then $L$ is the two-bridge link $\B(\alpha, \beta)$, \cite[Corollary 4.12]{hodgson-rubinstein:iaiols}.  

\begin{thm}\label{thm:ratlreps}
Let $L$ be an unoriented two-bridge link considered up to homeomorphism. \begin{enumerate}
\item No two-bridge link admits four equivalence classes of $3$--braid representatives.
\item $L$ admits three equivalence classes of $3$--braid representatives only if $L$ is equivalent to $\B(4,1)$.
\item $L$ admits two equivalence classes of $3$--braid representatives only if $L$ is equivalent to $\B(\alpha,1)$ and $\alpha \neq 0$
\item $L$ admits exactly one $3$--braid representative if $L$ is equivalent to either $\B(0,1)$ or $\B(\alpha, \beta)$ where $0<\beta<\alpha$ and either
\begin{itemize}
\item $\alpha=2pq+p+q$ and $\beta=2q+1$ for some integers $p,q>1$, or
\item $\alpha=2pq+p+q+1$ and $\beta=2q+1$ for some integers $p,q>0$.
\end{itemize}
\item $L$ admits no $3$--braid representatives otherwise.
\end{enumerate}
\end{thm}

The proof of Theorem~\ref{thm:ratlreps} is given in Section~\ref{sec:counting}.  Theorem~\ref{thm:generalgofks}, Theorem~\ref{thm:ratlreps}, and Corollary~4.12 of \cite{hodgson-rubinstein:iaiols} together imply the following corollary.

\begin{cor}\label{cor:countgofkils}
The lens space $L(\alpha, \beta)$ contains exactly zero, one, two, or three GOF-knots if and only if the two-bridge link $\B(\alpha, \beta)$ admits exactly zero, one, two, or three (respectively) equivalence classes of $3$--braid representatives as described in Theorem~\ref{thm:ratlreps}.  No lens space contains four GOF-knots.
\end{cor} 

\begin{remark}\label{remark:alternatetoH-R}
Since a link that may be represented as a closed $3$--braid has bridge number at most $3$, the full strength of \cite[Corollary 4.12]{hodgson-rubinstein:iaiols} is not necessary to obtain Corollary~\ref{cor:countgofkils}.  Birman and Hilden show that the double cover of $S^3$ branched over a link of bridge number $b \leq 3$ is a $3$--manifold of Heegaard genus $b-1$, \cite[Theorem 5]{birman-hilden:hsobcoS3}.  (As they remark, Viro independently proves this too, \cite{viro:ltsbcab}.)  Since, with the exception of $S^3$, lens spaces are the $3$--manifolds of Heegaard genus $1$, this implies that the only links with representations as closed $3$--braids that have a lens space as their double branched cover are two-bridge links and the unknot.
\end{remark}

%%%%%%%%%%%%%%%%%%
\section{Counting representations of two-bridge links as closed $3$-braids} \label{sec:counting}
Fortunately, most of the hard work for proving Theorem~\ref{thm:ratlreps} has been done.  Murasugi~\cite[Proposition 7.2]{murasugi:otbioal} and later Stoimenow~\cite[Corollary 8]{stoimenow:tspoc3b} determine which oriented two-bridge links have representations as closed $3$-braids.  The Classification Theorem of Birman and Menasco \cite{birman-menasco:slvcbIII} then permits us to count the number of braid axes representing an oriented two-bridge link as a closed $3$-braid that are not isotopic in the complement of the two-bridge link.  In Lemma~\ref{lem:distinctbraids} we show when these braid axes paired with the link are equivalent by a homeomorphism of $S^3$.   Theorem~\ref{thm:ratlreps} is then proved by determing which orientations of which two-bridge links admit closed $3$--braid representations. 

Let ${\bf b}(L)$ denote the braid index of the link $L$.
\begin{prop}[Proposition 7.2, \cite{murasugi:otbioal}] \label{prop:murasugi}
Let $L$ be a two-bridge link of type $\B(\alpha, \beta)$, where $0<\beta<\alpha$ and $\beta$ is odd.  Then
\begin{enumerate}
\item ${\bf b}(L) = 2$ iff $\beta = 1$.
\item ${\bf b}(L) = 3$ iff either
\begin{enumerate}
\item \label{item:one} for some $p,q >1$, $\alpha = 2pq+p+q$ and $\beta = 2q+1$, or
% the above line is a rewrite of the following
% \item for some $p,q>0$, $\alpha = 2pq+3p+3q+4$ and $\beta = 2p+3$, or
\item \label{item:two} for some $q>0$, $\alpha =2pq+p+q+1$ and $\beta = 2q+1$.
\end{enumerate}
\end{enumerate}
\end{prop} 

Since $\alpha$ is chosen to be positive and greater than $\beta$, the condition on type~\eqref{item:two} that $q>0$ forces $p>0$.

\begin{cor}[Corollary 8, \cite{stoimenow:tspoc3b}]\label{cor:stoimenow}
 If $L$ is a two-bridge link of braid index 3, then $L$ has Conway notation $(p,1,1,q)$ or $(p,2,q)$ for some $p,q>0$.
 \end{cor} 

\begin{figure}
\centering
\input{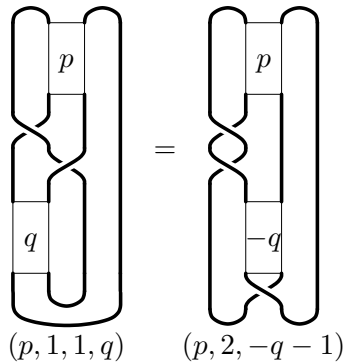}
\caption{The two bridge link $(p,1,1,q)$ is equivalent to $(p,2, -q-1)$.}
\label{fig:twobridgeequiv}
\end{figure}

\begin{remark}\label{remark:conwaynotation}
The link with Conway notation $(p,2,q)$ corresponds to type~\eqref{item:one} in Proposition~\ref{prop:murasugi}.
The link $(p,1,1,q)$ corresponds to type~\eqref{item:two} and is equivalent to the link $(p, 2, -q-1)$.  See Figure~\ref{fig:twobridgeequiv}.  Observe then that up to mirror equivalence the links $(p, 2, q)$ for any $p \in \Z_+$ and $q \in \Z$ contain all oriented two-bridge links of braid index at most $3$ and that every such link has braid index at most $3$.
\end{remark}

\begin{figure}
\centering
\input{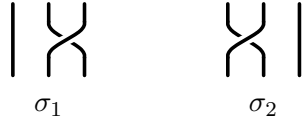}
\caption{The standard generators $\sigma_1$ and $\sigma_2$ for the three strand braid group $B_3$.}
\label{fig:standardgenerators}
\end{figure}

Let $\sigma_1$ and $\sigma_2$ be the standard generators of the $3$--braid group as depicted in Figure~\ref{fig:standardgenerators}.

\begin{thm}[The Classification Theorem, \cite{birman-menasco:slvcbIII}] \label{thm:classification}
 An oriented link $L$ which is represented by a closed $3$--braid admits a unique conjugacy class of $3$--braid representatives, with the following exceptions:
\begin{enumerate}
\item $L$ is the unknot, which has three conjugacy classes of $3$--braid representatives, namely the classes of $\sigma_1 \sigma_2$, $\sigma_1^{-1} \sigma_2^{-1}$ and $\sigma_1 \sigma_2^{-1}$.
\item $L$ is a type $(2, k)$ torus link, $k \neq \pm 1$, which has two conjugacy classes of $3$--braid representatives, namely the classes of $\sigma_1^k \sigma_2$ and $\sigma_1^k \sigma_2^{-1}$.
\item $L$ is one of a special class of links of braid index 3 which have $3$--braid representatives which admit ``braid-preserving flypes''.  These links have at most two conjugacy classes of $3$--braid representatives, namely the classes of $\sigma_1^p \sigma_2^r \sigma_1^q \sigma_2^\delta$ and $\sigma_1^p \sigma_2^\delta \sigma_1^q \sigma_2^r$, where $p$, $q$, $r$ are distinct integers having absolute value at least $2$ and where $\delta= \pm 1$.
\end{enumerate}
\end{thm}

\begin{figure}
\centering
\input{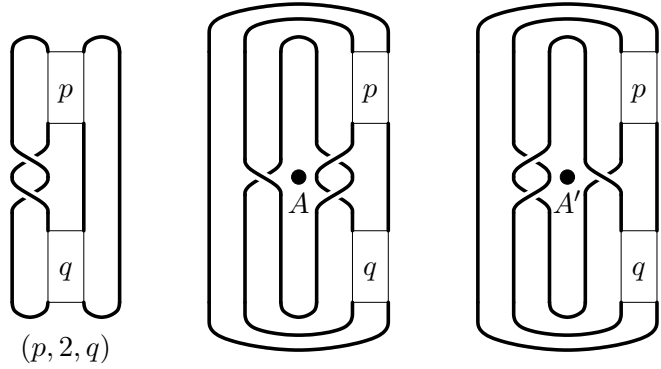}
\caption{The two bridge link $(p,2,q)$ and its two typically distinct conjugacy classes of closed $3$--braid representatives.}
\label{fig:p2qknot}
\end{figure}

\begin{remark}\label{remark:braidindex3reps}
As Figure~\ref{fig:p2qknot} shows, the two-bridge links with braid index $3$ (and hence links with Conway notation $(p, 2, q)$ up to mirror equivalence) belong to the third type of links in Theorem~\ref{thm:classification} where $r=2$ and $\delta = -1$.  Figure~\ref{fig:p2qknot} also indicates two braid axes $A$ and $A'$ for the two-bridge links of braid index $3$.  By Theorem~\ref{thm:classification} the braid axes $A$ and $A'$ are not isotopic in the complement of the two-bridge link if and only if $p, q \in \Z-\{-1, 0, 1, 2\}$ and $p \neq q$.  As shown in the proof of Lemma~\ref{lem:distinctbraids} below however, there is an involution of the two-bridge link that exchanges these two axes.
\end{remark}

\begin{lemma}\label{lem:distinctbraids}
An unoriented link $L$ which may be represented by a closed $3$--braid admits at most one equivalence class of braid axes giving $3$--braid representatives for a given orientation of $L$ and its reverse, with the following exception:  $L$ or its mirror is a type $(2,k)$ torus link with $k >0$, which has two equivalence classes of $3$--braid axes corresponding to the conjugacy classes of $\sigma_1^k \sigma_2$ and $\sigma_1^k \sigma_2^{-1}$ when coherently oriented .

\end{lemma}

\begin{figure}
\centering
\input{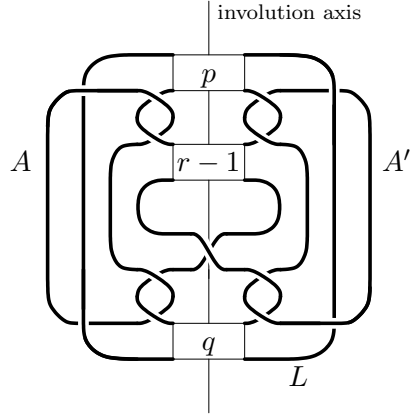}
\caption{The link $L$ of the third type in Theorem~\ref{thm:classification} (shown here with $\delta=-1$) admits an involution that exchanges its two braid axes $A$ and $A'$.}
\label{fig:involution}
\end{figure}

\begin{proof}
This lemma is perhaps suggested in \cite{birman-menasco:slvcbIII}.  
We only need consider the oriented links $L$ with at least two conjugacy classes of $3$--braid representatives as described in Theorem~\ref{thm:classification}.  

If $L$ is the unknot, the braid axes $A$ and $\bar{A}$ that correspond to the conjugacy classes of $\sigma_1 \sigma_2$ and $\sigma_1^{-1} \sigma_2^{-1}$ respectively are equivalent by an orientation reversing homeomorphism of $S^3$.  They are not equivalent, however, to the braid axis $A'$ that corresponds to the conjugacy class of $\sigma_1 \sigma_2^{-1}$ since the absolute values of the exponent sums of the braid words $\sigma_1 \sigma_2$ and $\sigma_1 \sigma_2^{-1}$ are not equal.  
Note that we may consider the unknot as a type $(2,1)$ torus link.

If $L$ is a type $(2, k)$ torus link with $k \neq \pm 1$, let $A$ and $A'$ be the two braid axes that correspond to the conjugacy classes of $\sigma_1^k \sigma_2$ and $\sigma_1^k \sigma_2^{-1}$.  If $k=0$, these axes are equivalent since there is an orientation reversing homeomorphism of $S^3$ taking $L$ to $L$ and $A$ to $A'$.  If $|k| \geq 2$ then these axes are not equivalent since the absolute values of the exponent sums of the braid words $\sigma_1^k \sigma_2$ and $\sigma_1^k \sigma_2^{-1}$ are not equal.  Since the $(2,-k)$ torus link is the mirror of the $(2,k)$ torus link, we may assume $k>0$.

If $L$ is a link that admits a ``braid-preserving flype,'' let $A$ and $A'$ be the braid axes that give (typically) distinct conjugacy classes of $3$--braid representatives.  The link $L$ and the axes $A$ and $A'$ are indicated in Figure~\ref{fig:involution} together with an axis of involution.  If $\iota$ is the involution about this axis, then $\iota(L)=L$ and $\iota(A)=A'$.   Hence the braid axes $A$ and $A'$ give equivalent $3$--braid representatives of $L$.  Note that $\iota$ reverses any orientation on $L$. 
\end{proof}

%%%%%%%%%%%%%%%%%%

\begin{proof}[Proof of Theorem~\ref{thm:ratlreps}]
By Lemma~\ref{lem:boundaxes}, an unoriented link $L$ has at most $4$ equivalence classes of braid axes that represent $L$ as a closed $3$--braid.  By Lemma~\ref{lem:orientedbraidreps} below, only the two-bridge link $\B(4,1)$ has two inequivalent orientations that each admit closed $3$--braid representatives.  Thus $\B(4,1)$ is the only two-bridge link that a priori could have more than two inequivalent closed $3$--braid representatives.

Theorem~\ref{thm:classification} implies that, as an oriented two-bridge link, $\B(4,1)$ has two closed $3$--braid representatives $\sigma_1^4 \sigma_2$ and $\sigma_1^4 \sigma_2^{-1}$ with axes $A_1$ and $A_2$ respectively.  By Lemma~\ref{lem:distinctbraids}, these two axes are inequivalent.  The other orientation, $\B(4,3)$, has Conway notation $(1,2,1)$ and a $3$--braid representative $\sigma_1 \sigma_2^2 \sigma_1 \sigma_2^{-1}$.  By Theorem~\ref{thm:classification} (and Lemma~\ref{lem:distinctbraids}), $\B(4,3)$ has just one closed $3$--braid representative with braid axis $A_3$.  The axes $A_1$ and $A_3$ are not equivalent since the absolute value of the exponent sums of their $3$--braid representatives are distinct.   
The axes $A_2$ and $A_3$ are not equivalent since $A_3$ cobounds an annulus with a component of $L$ whereas SnapPea \cite{weeks:snappea} reports $L \cup A_2$ as a hyperbolic link.
Thus the unoriented link $\B(4,1)$ admits a total of three equivalence classes of closed $3$--braid axes.  

Every unoriented two-bridge link of braid index $2$ is equivalent to $\B(\alpha, 1)$ as noted in Proposition~\ref{prop:murasugi}.  By Lemma~\ref{lem:orientedbraidreps} below and Lemma~\ref{lem:distinctbraids}, if $\alpha \neq 0$ or $4$ then such links have exactly two equivalence classes of braid axes giving closed $3$--braid representatives.  If $\alpha = 0$ then Lemma~\ref{lem:distinctbraids} implies that the link has just one equivalence class of braid axes representing $L$ as a closed $3$--braid.

By Proposition~\ref{prop:murasugi} a two-bridge link of braid index $3$ is equivalent to $\B(\alpha, \beta)$ with $0 < \beta < \alpha$  if and only if $\alpha$ and $\beta$ satisfy either \eqref{item:one} or \eqref{item:two} of the proposition.  By Lemma~\ref{lem:orientedbraidreps} and Lemma~\ref{lem:distinctbraids} these links have just one equivalence class of braid axes representing $L$ as a closed $3$--braid.

By Proposition~\ref{prop:murasugi} a two-bridge link has no closed $3$--braid representatives if it is not equivalent to some $\B(\alpha,\beta)$ where either $\beta=1$, \eqref{item:one} is satisfied, or \eqref{item:two} is satisfied.  Thus such a link has no equivalence classes of braid axes representing $L$ as a closed $3$--braid.
\end{proof}

\begin{remark} \label{rem:twobridgeequiv}
Recall that two oriented two-bridge links $\B(\alpha_1, \beta_1)$ and $\B(\alpha_2, \beta_2)$ are ambient isotopy  equivalent as oriented links if and only if
\[ \alpha_1 = \alpha_2  \mbox{   and   }  \beta_1^{\pm1} \equiv \beta_2 \mod 2 \alpha_1 \]
whereas the unoriented two-bridge links are equivalent if and only if the second condition is taken simply mod $\alpha_1$.

Since oriented two-bridge links are invertible, there is no distinction between the orientations of a two-bridge knot.  However, switching the orientation on one component of the oriented link $\B(\alpha, \beta)$ yields the link $\B(\alpha, \beta-\alpha)$ which is mirror equivalent to $\B(\alpha, \alpha-\beta)$.  Thus, up to mirror equivalence, every oriented two-bridge link except the two component unlink $\B(0,1)$ and the unknot $\B(1,1)$ is represented by $\B(\alpha, \beta)$ for some $0 < \beta < \alpha$, with $\beta$ odd; see, e.g., \cite{bz:knots}.
\end{remark}

\begin{remark}\label{rem:murasugiequiv}
If $\beta$ and $\beta'$ are odd integers with $0<\beta<\alpha$ such that $\beta \beta' \equiv 1 \mod 2 \alpha$, then for positive integers $p$ and $q$
\begin{itemize}
\item if $\alpha = 2pq + p + q$  and $\beta = 2p+1$, then $\beta' = 2q+1$, and
\item if $\alpha = 2pq + p + q + 1$  and $\beta = 2p+1$, then $\beta' =-(2q+1)$.
\end{itemize}
Therefore, up to mirror equivalence, the oriented two-bridge link $\B(2pq+p+q+\delta, 2p+1)$ is equivalent to $\B(2pq+p+q+\delta, 2q+1)$ where $\delta \in \{0,1\}$.
\end{remark}

\begin{lemma}\label{lem:orientedbraidreps}
Among two-bridge links, only the link $\B(4,1)$ has two distinct orientations which each admit a closed $3$--braid representation.
\end{lemma}

\begin{proof}
Assume the unoriented two-bridge link $L$ has two distinct orientations and a closed $3$--braid representative for each.  Then, as noted in Remark~\ref{rem:twobridgeequiv}, $L$ is necessarily a link of two components.  Thus up to mirror equivalence, the two orientations of $L$ may be denoted as $\B(\alpha, \beta)$ and $\B(\alpha, \alpha-\beta)$ where $0 < \beta < \alpha$.  Since $\alpha$ is necessarily even, both $\beta$ and $\alpha-\beta$ are odd.  

\noindent {\em Case 1.}  $\beta = 1$ or $\alpha-\beta = 1$

We may assume $\beta = 1$.  If $\alpha = 2$ then the two orientations on $\B(2,1)$ are mirror equivalent.  Hence we may further assume $\alpha >2$.

Theorem~\ref{thm:classification} shows that $\B(\alpha, 1)$ has a closed $3$--braid representative.  %, in fact exactly two inequivalent $3$--braid representatives, $\sigma_1^\alpha \sigma_2$ and $\sigma_1^\alpha \sigma_2^{-1}$.  
Since $\alpha >2$, a $3$--braid representative of $\B(\alpha, \alpha-1)$ must be of the second type in Proposition~\ref{prop:murasugi}.  Therefore, in accordance with Remark~\ref{rem:murasugiequiv}, we only need check if 
\[\alpha = 2pq + p + q + \delta \mbox{   and   } \alpha-1 = 2p+1 \]
for some integers $p,q>0$ and $\delta \in \{0,1\}$.  
It follows that $\alpha = 2p+2$ and hence $p=(2-q-\delta)/(2q-1)$.  The only valid solution is $p=1=q$ with $\delta = 0$.  Thus $\alpha = 4$.  Because $1 \cdot 3 \equiv 3 \mod 2 \cdot 4$, the two oriented links $\B(4,1)$ and $\B(4,3)$ are not mirror equivalent.  Therefore each orientation of the unoriented two-bridge link $L = b(4,1)$ has a $3$--braid representative.

\noindent {\em Case 2.}  $\beta >1$ and $\alpha-\beta >1$

Any $3$--braid representative of $L$ must be of the second type in Proposition~\ref{prop:murasugi}.  Therefore, again in accordance with Remark~\ref{rem:murasugiequiv}, we only need check if 
\[\alpha = 2pq + p + q + \delta \mbox{   and   } \beta = 2p+1 \]
and
\[\alpha = 2rs + r + s + \epsilon \mbox{   and   } \alpha-\beta = 2r+1 \]
for some integers $p, q, r, s>0$ and $\delta, \epsilon \in \{0,1\}$.
Eliminating $\beta$, we have the three equations 
\begin{align}
\alpha &= 2pq + p + q + \delta = (2p+1)q + p + \delta \label{eq:first}\\
\alpha &= 2rs + r + s + \epsilon = (2r+1)s + r + \epsilon \label{eq:second}\\
\alpha &= (2p+1) + (2r+1). \label{eq:third}
\end{align}
Combining Equation~\eqref{eq:third} with each \eqref{eq:first} and \eqref{eq:second} we obtain
\begin{align}
(2r+1) &= (2p+1)(q-1) + p + \delta  \label{eq:fourth}\\
(2p+1) &= (2r+1)(s-1) + r + \epsilon. \label{eq:fifth}
\end{align}
By examining Equation~\eqref{eq:fourth}, if $q = 1$ then $p>r$ and if $q >1$ then $p<r$.  Similarly Equation~\eqref{eq:fifth} implies that if $s = 1$ then $r>p$ and if $s>1$ then $p>r$.  Hence either $q=1$ and $s>1$ (in which case $p>r$) or $q>1$ and $s=1$ (in which case $r>p$).  These two cases are symmetric.

Assume $q=1$ and $s>1$.  Then Equation~\eqref{eq:fourth} gives
\begin{align} 2r+1 = p + \delta. \label{eq:sixth} \end{align}
Substituting this into Equation~\eqref{eq:fifth} yields
\begin{align*}
2p+1 &= (p+\delta)(s-1) + r + \epsilon \\
         &= (s-1)p + (s-1) \delta + r + \epsilon 
\end{align*}
and thus
 \begin{align} 
 (3-s)p + 1= (s-1)\delta + r + \epsilon. \label{eq:seventh} 
 \end{align}
 Since the right hand side is necessarily positive, $s=2$ or $s=3$.
 
 If $s=2$, then 
 \begin{align*} 
 p+1 &= r+\delta + \epsilon \mbox{  by Eq.~\eqref{eq:seventh}} \\
 2r+1-\delta +1 &= r+\delta + \epsilon \mbox{  by Eq.~\eqref{eq:sixth}} \\
r &= 2\delta + \epsilon -2
\end{align*}
Since $r>0$, $\delta=\epsilon=1$, $r=1$, and $p=2$.  Thus $\alpha = 8$, $\beta =5$, and $\alpha-\beta=3$.  Because $3 \cdot 5 \equiv -1 \mod 2 \cdot 8$, the oriented links $\B(8,3)$ and $\B(8,5)$ are mirror equivalent.  %Note that having $r=q$ and $p=s$ evidences this.

If $s=3$, then by Equation~\eqref{eq:seventh},
\[1 = 2 \delta + r + \epsilon.\]
Since $r>0$, $\delta=\epsilon=0$, $r=1$, and $p=3$.  Thus $\alpha=10$, $\beta=7$, and $\alpha-\beta=3$.  Because $3 \cdot 7 \equiv 1 \mod 2 \cdot 10$, the oriented links $\B(10,3)$ and $\B(10,7)$ are equivalent.
\end{proof}

%%%%%%%%%%%%%%%%%%
\section{Remarks}

\begin{remark}
One may obtain explicit pictures of the fiber surface of these GOF-knots in lens spaces   like those in \cite{morimoto:gofkils} by carrying a disk that both is bounded by the braid axis and intersects the $3$--braid minimally through the sequence of steps done to obtain a presentation of a lens space as surgery on the unknot from its corresponding two-bridge link.
\end{remark}

\begin{remark}
As Morimoto notes in \cite[Remark 1]{morimoto:gofkils}, his knot $K_2$ in the lens space $L(5,1)$ has two meridians.  A $+1$ surgery on $K_2$ produces the manifold $L(5,4)$ which is equivalent to $L(5,1)$ by an orientation reversing homeomorphism.  

We observe this in the context of closed $3$--braids as follows.  The knot $K_2$ in $L(5,1)$ is the lift of the braid axis in the double branched cover of the closure of the braid $\sigma_1^5 \sigma_2$, the two-bridge knot $\B(5,1)$.   As noted in Remark~\ref{remark:axissurgery}, $+1$ surgery on $K_2$ corresponds to inserting two full left-handed twists (i.e.\ $((\sigma_1 \sigma_2)^3)^{-2}$) into the braid.  
 Therefore $+1$ surgery on $K_2$ corresponds to the double branched cover of the closure of the braid
\[
(\sigma_1^5 \sigma_2) ((\sigma_1 \sigma_2)^3)^{-2} 
= (\sigma_1^5 \sigma_2)(\sigma_2^{-1} \sigma_1^{-3} \sigma_2^{-2} \sigma_1^{-1} \sigma_2^{-3} \sigma_1^{-2})
\equiv \sigma_2^{-5} \sigma_1^{-1} 
\equiv \sigma_1^{-5} \sigma_2^{-1}
\]
where $\equiv$ denotes equivalence up to conjugation.  
 The closure of $\sigma_1^{-5} \sigma_2^{-1}$ may be recognized as the two-bridge knot $\B(5,4)$ and the mirror of the closure of $\sigma_1^5 \sigma_2$.
%\end{remark}

%\begin{remark}
Because $1/n$ surgery on a GOF-knot confers a GOF-knot in the surgered manifold, the above example appears to be the only situation in which $1/n$ surgery on a GOF-knot in a lens space yields a homeomorphic lens space.  Clearly this example may be generalized to obtain GOF-knots in manifolds other than lens spaces that admit a $1/n$ surgery yielding a homeomorphic manifold with the opposite orientation.  See \cite{bhw:csok} for more on this sort of phenomenon.
\end{remark}

\begin{remark}
As Morimoto shows in \cite{morimoto:gofkils} and is further observed in Corollary~\ref{cor:countgofkils}, there are lens spaces which contain no GOF-knots.  Nevertheless, there are knots in lens spaces representing a non-trivial element of homology whose exteriors are once-punctured torus bundles.  For instance, since $-17/3$ surgery on the right-handed trefoil yields $L(17,5)$, even though the core of the surgered solid torus is not a GOF-knot in $L(17,5)$, its exterior is a once-punctured torus bundle.  Indeed Corollary~\ref{cor:countgofkils} shows that $L(17,5)$ contains no GOF-knots.
\begin{question}
Must a lens space contain a knot whose exterior is a once-punctured torus bundle?  If not, which do?
%What lens spaces have no GOF-knots but contain a knot whose exterior is a once-punctured torus bundle?
\end{question}
\end{remark}

\begin{remark}
Via double branched covers, one obtains a genus $g$ fibered knot from the braid axis of a closed braid of braid index $2g+1$.  However, not all genus $g$ fibered knots arise in this manner if $g>1$.  Again considering that there are lens spaces that contain no GOF-knots, we ask the following question.
\begin{question}
What is the minimal genus among (null homologous) fibered knots in a given lens space $L(\alpha,\beta)$?
\end{question} 
Murasugi shows that the braid index of a two-bridge link may be arbitrarily large, \cite[Theorem B]{murasugi:otbioal}.   Perhaps then it is not too foolish to conjecture that there exist lens spaces whose minimal genus (null homologous) fibered knot has arbitrarily large genus.
\end{remark}

The author would like to thank Will Kazez for useful conversations and Kanji Morimoto for his provoking article.

\bibliography{GOFKBiblio}
\bibliographystyle{hamsplain}

\end{document}